\documentstyle[12pt]{article}
\def\N{I\!\!N}
\def\bbbn{I\!\!N}
\def\bbbr{\rm I\!R}
\def\bbbc{C\!\!\!\!I}
\def\bbbq{Q\!\!\!\!I}
\def\bbbz{Z\!\!\!Z}
\def\Z{Z\!\!\!Z}
\def\Tr{\rm Tr\ }
\def\GL{\rm GL }
\def\U{\hbox{U}}

\def\Fix{\rm {\hbox{Fix}}}

\def\cR{{\cal R}}
\def\cS{{\cal S}}

\newtheorem{theorem}{Theorem}
\newtheorem{lemma}{Lemma}
\newtheorem{corollary}{Corollary}
\newtheorem{proposition}{Proposition}
\newtheorem{definition}{Definition}
\newtheorem{example}{Example}

\newcommand{\ind}{{\rm Index\ }}
\newcommand{\fix}{{\rm Fix\ }}
\newcommand{\M}{{\rm M}}
\newcommand{\mod}{{\rm \; mod\ }}
\newcommand{\proof}{{\it Proof. }}
\newcommand{\pr}{{\rm pr\ }}
\newcommand{\rank}{{\rm Rank\ }}

\newcommand{\spec}{{\rm Spec\ }}
\newcommand{\square}{\hfill$\Box$}
\newcommand{\tr}{{\rm Tr\ }}

\date{}
\title{Trace Formulae, Zeta Functions, Congruences and Reidemeister 
Torsion in Nielsen Theory}
 \author{Alexander Fel'shtyn\thanks{Part of this work was conducted during
 authors' stay 
 in Sonderforschungsbereich 170 ``Geometrie und Analysis'', 
 Mathematisches Institut der Georg August Universit\"at zu 
 G\"ottingen.} \and Richard Hill}

\begin{document}
\bibliography{ref}
\bibliographystyle{plain}
\maketitle

\bigskip
\begin{abstract}
In this paper we prove trace formulae for the Reidemeister number of 
 a group endomorphism.
This result implies the rationality of the Reidemeister zeta 
 function in the following cases:
 the group is a direct product of a finite group and a finitely generated
Abelian group;
 the group is finitely generated, nilpotent and torsion free. 
%We continue to study analytical properties of the Nielsen zeta function. 
We connect the Reidemeister zeta function of an endomorphism of a 
 direct product of a finite group and a finitely generated free Abelian 
 group with the Lefschetz zeta function of the induced map on the unitary 
 dual of the group.
As a consequence we obtain a relation
 between a special value of the Reidemeister zeta function
 and a certain Reidemeister torsion.
We also prove congruences for Reidemeister numbers of iterates of an
endomorphism
 of a direct product of a finite group and a finitely generated free Abelian
group
 which are the same as those found by Dold for Lefschetz numbers.
\end{abstract}

\setcounter{section}{-1}

\section{Introduction}
We assume everywhere $X$ to be a connected, compact polyhedron and 
 $f:X\rightarrow X$ to be a continuous map.
Taking a dynamical point of view, we shall be interested in the iterates of
$f$.
In the theory of discrete dynamical systems the following zeta functions 
 have been studied: the Artin-Mazur zeta function
$$
 \zeta_f(z) := \exp\left(\sum_{n=1}^\infty \frac{F(f^n)}{n} z^n \right),
$$
 where $F(f^n)$ is the number of isolated fixed points of $f^n$;
 the Ruelle zeta function \cite{ruel}
$$
{\zeta_f}^g(z) := \exp\left(\sum_{n=1}^\infty \frac{ z^n 
}{n}\sum_{x\in \Fix (f^n)} \prod_{k=0}^{n-1}g(f^k(x)) \right),
$$
 where $g:X\to \bbbc $ is a weight function (if $g$ is identically 1 
 then this is simply the Artin-Mazur function); 
 the Lefschetz zeta function
$$
 L_f(z) := \exp\left(\sum_{n=1}^\infty \frac{L(f^n)}{n} z^n \right),
$$
 where
$$
 L(f^n) := \sum_{k=0}^{\dim X} (-1)^k \tr\Big[f_{*k}^n:H_k(X;\bbbq)\to 
 H_k(X;\bbbq)\Big]
$$
 are the Lefschetz numbers of the iterates of $f$;
 reduced modulo 2 Artin-Mazur and Lefschetz zeta functions \cite{fran};
 twisted Artin-Mazur and Lefschetz zeta functions \cite{fri},
 which have coefficients in the group ring $\bbbz H$ of an Abelian group $H$.

The above zeta functions are analogous to the Hasse-Weil 
 zeta function of an algebraic variety over a finite field \cite{weil}.
Like the Hasse-Weil zeta function, the Lefschetz function
 is always a rational function of $z$, and is given by the formula:
$$
 L_f(z) = \prod_{k=0}^{\dim X} \det\big(I - f_{*k}.z\big)^{(-1)^{k+1}}.
$$
This immediately follows from the trace formula for the Lefschetz numbers 
 of the iterates of $f$.
The Artin-Mazur zeta function has a positive radius of convergence
 for a dense set in the space of smooth self-maps
 of a compact smooth manifold \cite{am}.
Manning proved the rationality of the Artin-Mazur zeta function
 for diffeomorphisms of a compact smooth manifold satisfying Smale's Axiom
A \cite{m}.
The value of knowing that a zeta function is rational is that it shows 
 that the infinite sequence of coefficients is closely interconnected, 
 and is given by the finite set of zeros and poles of zeta function. 

The Artin-Mazur zeta function and its modifications count periodic 
 points of a map geometrically, and the Lefschetz's type zeta functions do 
 the same thing algebraically (with weight defined by index theory).
A third way of counting periodic points is to use Nielsen theory.
Let $p:\tilde{X}\rightarrow X$ be the universal covering of $X$ and 
 $\tilde{f}:\tilde{X}\rightarrow \tilde{X}$ a lifting of $f$, ie. 
 $p\circ\tilde{f}=f\circ p$.
Two liftings $\tilde{f}$ and $\tilde{f}^\prime$ are said to be {\it
conjugate}
 if there is a $\gamma\in\Gamma\cong\pi_1(X)$ such that 
 $\tilde{f}^\prime = \gamma\circ\tilde{f}\circ\gamma^{-1}$.
The subset $p(\Fix (\tilde{f}))\subset \Fix (f)$ is called
 {\it the fixed point class of $f$ determined by the lifting class
$[\tilde{f}]$}.
A fixed point class is said to be {\it essential} if its index is non-zero.

The number of lifting classes of $f$ (and hence the number of fixed 
 point classes, empty or not) is called the {\it Reidemeister number} 
 of $f$, denoted $R(f)$.
This is a positive integer or infinity.
The number of essential fixed point classes is called
 the {\it Nielsen number} of $f$, denoted by $N(f)$.
The Nielsen number is always finite.
Both Nielsen and Reidemeister numbers are homotopy type invariants.
In the category of compact, connected polyhedra the Nielsen number of 
 a map is a lower bound  for the least number of fixed points of maps in 
 the homotopy class of $f$.

Let $G$ be a group and $\phi:G\rightarrow G$ an endomorphism.
Two elements $\alpha,\alpha^\prime\in G$ are said to be $\phi-conjugate$ 
 if there exists $\gamma \in G$ such that $\alpha^\prime=\gamma \cdot 
 \alpha \cdot \phi(\gamma)^{-1}$.
The number of $\phi$-conjugacy classes is called
 the {\it Reidemeister number} of $\phi$, denoted by $R(\phi)$.
We shall also write $\cR(\phi)$ for the set of $\phi$-conjugacy classes
 of elements of $G$.
One easily shows that if $\psi$ is an inner automorphism of $G$ then
 $R(\phi\circ\psi)=R(\psi\circ\phi)=R(\phi)$.
The group-theoretical and topological Reidemeister numbers are related as
follows.
If $G$ is the fundamental group of $X$ and $\phi$ is the endomorphism
 of $G$ induced by $f$ (this is only defined modulo inner 
 automorphisms) then $R(\phi)=R(f)$.

If we consider the iterates of $f$ and $\phi$, we may define several 
 zeta functions associated with Nielsen fixed point theory (see 
 \cite{f4,f5,fp}).
We assume throughout this article that $R(f^n)<\infty$ and
$R(\phi^n)<\infty$
 for all $n>0$.
The Reidemeister zeta functions of $f$ and $\phi$ and the Nielsen 
zeta function of $f$ are defined as power series:
\begin{eqnarray*}
 R_\phi(z) & := & \exp\left(\sum_{n=1}^\infty \frac{R(\phi^n)}{n} z^n
\right), \\
 R_f(z) & := & \exp\left(\sum_{n=1}^\infty  \frac{R(f^n)}{n} z^n \right),
\\
 N_f(z) & := &  \exp\left(\sum_{n=1}^\infty \frac{N(f^n)}{n} z^n \right). 
\end{eqnarray*}
These functions are homotopy invariants.
The Nielsen zeta function $N_f(z)$ has a positive radius of convergence
which has a 
 sharp estimate in terms of the topological entropy of the map $f$ \cite{fp}.
In section 2 we give another proof of positivity of the 
 radius of convergence.
We also obtain an exact algebraic lower bound for the radius using 
 the Reidemeister trace formula for generalized Lefschetz numbers.

We begin the article by proving in section 1 trace formulae for 
 Reidemeister numbers in the following cases:
 $G$ is finite;
 $G$ is a direct product of a finite group and a finitely generated free 
Abelian group;
 $G$ is a finitely generated torsion free nilpotent group.
By this we mean that there are finite dimensional complex vector 
 spaces $H^{i}$ for $i=0,\ldots N$ and linear maps $T_{i}:H^{i}\to H^{i}$
 such that for every natural number $n$, the Reidemeister number of
$\phi^{n}$
 is equal to $\sum_{i=0}^{N}(-1)^{n}\Tr(T_{i}^{n})$.
Thus the definition of the Lefschetz number is a trace formula.
These results had previously been known only for the finitely 
 generated free Abelian groups \cite{f5}, although the case of finite 
 groups is implicit in \cite{fh1,fh2}

The trace formula for the Reidemeister numbers implies the rationality of 
 the Reidemeister zeta function.
In section 2 we also prove functional equations
 for the Reidemeister zeta function of group endomorphisms in the 
 above cases. 

In section 2 we  also study analytic properties of the Nielsen zeta function.
We give a new proof of the positivity of the radius of convergence of
$N_f(z)$,
 and give an exact algebraic lower bound for the radius.
The radius of convergence of $N_f(z)$ is directly related
 to the asymptotic Nielsen number, which is defined to be the growth rate
 of the sequence $N(f^n)$ of the Nielsen numbers of the iterates of $f$. 
It is worth mentioning that in all known cases
 the Nielsen zeta function is a nice function.
By this we mean that it is a product of 
 an exponential of a polynomial with a function some power of which is
rational.

In section 1 we prove arithmetical congruences :
$$
 \sum_{d\mid n} \mu(d)\cdot R(\phi^{n/d}) \equiv 0 \mod n,
$$
 for the Reidemeister numbers of the iterations of a group 
 endomorphism $\phi$ of a direct product of a finite group and a finitely
generated free 
 Abelian group and corresponding congruences
$$
 \sum_{d\mid n} \mu(d)\cdot R(f^{n/d}) \equiv 0 \mod n
$$
 for the Reidemeister numbers of a continuous map $f$.
These congruences are the same as those found by Dold \cite{do}
 for Lefschetz numbers. 

Recently a connection between the Lefschetz type dynamical zeta 
 functions and the Reidemeister torsion was established by D. Fried
\cite{fri1,fri2}.
The work of Milnor \cite{mi} was the first indication that such a
connection exists. 
In section 3 we establish a connection between Reidemeister 
 torsion and the Reidemeister zeta function.
We obtain an expression for the Reidemeister torsion of the mapping torus
 of the dual map of an endomorphism of a direct product of a finite group
 and a finitely generated free Abelian group,
 in terms of a special value of the Reidemeister zeta function of the
endomorphism.
The result is obtained by expressing the Reidemeister zeta function
 in terms of the Lefschetz zeta function of the dual map,
 and then applying the theorem of D. Fried.
This means  that the Reidemeister torsion counts the periodic point
classes
 of the map $f$.
These results had previously been known for the finitely generated Abelian
groups
 and finite groups\cite{fh2}.

We would like to thank J. Eichhorn , D. Fried, M. L. Gromov and B. B. 
Venkov for valuable conversations and comments. We would like to 
thank Institut des Hautes Etudes Scientifiques in Bures-sur-Yvette 
and the Max Planck Institut f\"ur Mathematik in Bonn for their kind 
hospitality and support.

\section[Group Endomorphisms]{Group Theoretical Reidemeister Numbers} 

Let $\phi$ be an endomorphism of a group $G$.
We shall write the group law of $G$ multiplicatively.
We shall write $\{g\}_{\phi}$ for the $\phi$-conjugacy class of an element
$g\in G$.
We shall write $<g>$ for the ordinary conjugacy class of $g$ in $G$.
The endomorphism $\phi$ maps conjugate elements to conjugate elements.
It therefore induces an endomorphism of the set of conjugacy classes
 in $G$.

\begin{theorem}[\cite{fh2}]
Let $G$ be a finite group and let $\phi:G\to G$ be an endomorphism. 
Then $R(\phi)$ is the number of ordinary conjugacy classes $<x>$ in 
$G$ such that
$$
 <\phi(x)>=<x>.
$$
\end{theorem}

We now rephrase this to give a trace formula.
Suppose $G$ is finite and let $W$ be the vector space of complex valued
 class functions on the group $G$.
A class function is a function which takes the same 
 value on every element of a usual conjugacy class.
The map $\phi$ induces a linear map $B:W\rightarrow W$ defined by
$$
 B(f):=f\circ\phi.
$$
We shall write $\hat G$ for the set of isomorphism classes
 of irreducible unitary representations of $G$.
There is a multivalued map $\hat\phi:\hat G\to \hat G$
 defined as follows.
Let $\rho:G\to\U(V)$ be an irreducible unitary representation.
Then we obtain a representation $\phi^{*}\rho$ of $G$ on $V$ by
 $\phi^{*}\rho(g)=\rho(\phi(g))$.
The representation $\phi^{*}\rho$ is a sum of irreducible pieces;
 we define $\hat\phi(\rho)$ to be this set of pieces.
If $\phi$ is bijective or $G$ is Abelian then $\hat\phi$ is single valued.

\begin{theorem}
\label{finitetrace}
Let $\phi:G\rightarrow G$ be an endomorphism of a finite group $G$. 
Then we have
\begin{equation}
 R(\phi)
 =
 \#\Fix(\hat\phi)
 =
 \Tr B
\end{equation}
\end{theorem}

In the case when $\phi$ is the identity map, this is Burnside's result 
 equating the number of irreducible representation of a finite group 
 with the number of conjugacy classes in the group.
\medskip

\proof
We shall calculate the trace of $B$ in two ways.
The characteristic functions of the conjugacy classes in $G$ form a basis
of $W$,
 and are mapped to one another by $B$
 (the map need not be a bijection).
Therefore the trace of $B$ is the number of elements 
 of this basis which are fixed by $B$.
By Theorem 1, this is equal to the Reidemeister number of $\phi$. 
Another basis of $W$, which is also mapped to itself by $B$ is the 
 set of traces of irreducible representations of $G$ (see \cite{la} chapter
XVIII). 
From this it follows that the trace of $B$ 
 is the number of irreducible representations $\rho$ of $G$ such that 
 $\rho$ has the same trace as $\hat\phi(\rho)$.
However, representations of finite groups 
 are characterized up to equivalence by their traces.
Therefore the trace of $B$ is equal to the number of fixed points of
$\hat\phi$.
\square
\medskip

The following result is known for free Abelian groups:

\begin{theorem}[\cite{fh2}]
\label{determinant}
Let $\phi:\Z^{k}\to\Z^{k}$ be a group endomorphism.
Then
$$
 R(\phi) = |\det(I-\phi)| = \#\Fix(\hat\phi).
$$
Furthermore the following trace formula holds:
\begin{equation}
 R(\phi) = (-1)^{r+p} \sum_{i=0}^k (-1)^i \tr (\wedge^i\phi). 
\end{equation}
 where $p$ the number of eigenvalues $\mu$ of $M$ such that $\mu <-1$,
 and $r$ the number of real eigenvalues whose absolute value is $>1$.
Here $\wedge^i$ denotes the exterior power.
\end{theorem}

Now let $F$ be a finite group and $k$ a natural number. 
We shall consider an endomorphism $\phi$ of the group $G = \bbbz^k \times
F$.
Our aim is to prove a trace formula for the Reidemeister number of 
 such an endomorphism.
The torsion elements of $G$ 
 are precisely the elements of the finite, normal subgroup $F$.
For this reason we have $\phi(F) \subset F$.
Let $\phi^{finite}: F\rightarrow F $
 be the restriction of $\phi$ to $F$,
 and let $\phi^{\infty}: G/F \rightarrow G/F$
 be the induced map on the quotient group. 
Let $\pr_{\bbbz^k}: G\rightarrow \bbbz^k$
 and $\pr_F :G\rightarrow F$ denote the projections onto $\bbbz^k$ and $F$.
Then the composition
$$
 \pr_{\bbbz^k}\circ\phi: \bbbz^k \rightarrow G\rightarrow \bbbz^k
$$
 is an endomorphism of $\bbbz^k$,
 which is given by some matrix $M \in \M_k(\bbbz)$.
We denote by $\psi : \bbbz^k \rightarrow F$
 the other component of the restriction of $\phi$ to $\bbbz^k$, i.e. 
$$
 \psi(v)=\pr_F(\phi(v)).
$$
We therefore have for any element $(v,f) \in G$
$$
 \phi(v,f)=(M\cdot v,\psi(v)\phi(f)).
$$

\begin{lemma}
Let $G$ be as above.
Two elements $g_1=(v_1,f_1)$ and $g_2=(v_2,f_2)$ of $G$.
 are $\phi$-conjugate if and only if
$$
 v_1\equiv v_2 \mod (1-M)\bbbz^k
$$
 and there is a $h\in F$ with
$$
 hf_1=f_2 \phi((1-M)^{-1}(v_2-v_1))\phi(h).
$$
\end{lemma}

\proof
Suppose $g_1$ and $g_2$ are $\phi$-conjugate.
Then there is a $g_3=(w,h) \in G$ with $g_3g_1= g_2\phi(g_3).$
Therefore
$$
 (w+v_1, hf_1)=(v_2+ M\cdot w, f_2\psi(w)\phi(h)).
$$
Comparing the first components we obtain $(1-M)\cdot w=v_2 -v_1$ 
 from which it follows that $v_1$ is congruent to $v_2$ modulo 
 $(1-M)\bbbz^r$.
Substituting $(1-M)^{-1}(v_2 -v_1)$ for $w$ in the 
 second component we obtain the second relation in the lemma.
The argument can easily be reversed to give the converse. 
\square

\begin{proposition}
\label{product}
In the notation described above,
 $R(\phi)= R(\phi^{finite}) \times R(\phi^{\infty})$.
\end{proposition}

\proof
We partition the set $\cR(\phi)$ of $\phi$-conjugacy classes of 
 elements of $G$ into smaller sets:
$$
 \cR(\phi)= \cup_{v\in \bbbz^k/(1-M)\bbbz^k} \cR(v),
$$
 where $\cR(v)$ is the set of $\phi$-conjugacy classes 
 $\{(w,f)\}_{\phi}$ for which $w$ is congruent to $v$
 modulo $(1-M)\bbbz^k$.
It follows from the previous lemma that this is a partition.
Now suppose $\{(w,f)\}_{\phi} \in \cR(v) $.
We will show that $\{(w,f)\}_{\phi}=\{(v,f^*)\}_{\phi}$ for some $f^* \in
F$.
This follows by setting $ f^*=f \psi((1-M)^{-1}(w-v))$
 and applying the previous lemma with $h=id$.
Therefore $\cR(v)$ is the set of $\phi$-conjugacy classes
$\{(v,f)\}_{\phi}$
 with $f \in F$.
From the previous lemma it follows that $(v,f_1)$ and $(v,f_2)$ are 
 $\phi$-conjugate iff there is a $h\in F$ with
$$
 hf_1=f_2\psi(0)\phi(h)=f_2\phi(h)
$$
This just means that $f_1$ and $f_2$ are $\phi^{finite}$-conjugate as 
 elements of $F$.
From this it follows that $\cR(v)$ has cardinality 
 $R(\phi^{finite})$.
Since this is independent of $v$, we have
$$
 R(\phi)
 =
 \sum_v R(\phi^{finite})
 =
 \mid \det(1-M)\mid \times R(\phi^{finite}).
$$

Now consider the map $\phi^{\infty}: G/F \rightarrow G/F $.
We have
$$
 \phi^{\infty}((v,F))=(M\cdot v, \psi(v)F)= (M\cdot v,F).
$$
From this it follows that $\phi^{\infty}$
 is equivalent to the map $M:\bbbz^k \rightarrow \bbbz^k$.
This implies
$$
 R(\phi^{\infty}) = R(M:\bbbz^k \rightarrow \bbbz^k).
$$
However by Theorem \ref{determinant} we have
 $R(M:\bbbz^k \rightarrow \bbbz^k)=\mid \det(1-M)\mid$.
Therefore $R(\phi)= R(\phi^{finite}) \times R(\phi ^{\infty})$,
 proving Proposition 1. 
\square
\medskip

\begin{proposition}
In the notation described above 
$$
 \#\fix(\hat\phi)
 =
 \#\fix(\hat\phi^{finite})\times 
 \#\fix(\hat\phi^\infty)
$$
\end{proposition}

\proof
Consider the dual $\hat G$. This is cartesian product of the duals 
of $ \bbbz^k$ and $F$:
$$
 \hat G=\hat \bbbz^k\times\hat F,   \>  \rho=\rho_1\otimes \rho_2
$$
 where $\rho_1$ is an irreducible representation of $\bbbz^k$ and 
 $\rho_2$ is an irreducible
 representation of $F$.
Since $\bbbz^r$ is Abelian, all of its irreducible  representations are
 1-dimensional, so $\rho(v)$ for $v\in \bbbz^k$ is always a scalar 
 matrix, and $\rho_2$ is the restriction of $\rho$ to $F$.
If $\hat\phi(\rho)=\rho$ then 
 there is a matrix $T$ such that 
$$
 \rho\circ\phi=T\cdot \rho \cdot T^{-1}.
$$
This implies 
$$
 \rho^{finite}\circ\phi^{finite}=T\cdot\rho^{finite} \cdot T^{-1},
$$
 so $\rho_2=\rho^{finite}$ is fixed by $\hat\phi^{finite}$.
For any fixed $\rho_2\in \cS(\phi^{finite})$, the set of $\rho_1$ 
 such that $\rho_1\otimes \rho_2$ is fixed by $\hat\phi$
 is the set of $\rho_1$ satisfying 
$$
 \rho_1(M\cdot v)\rho_2(\psi(v))= T\cdot \rho_1(v) \cdot T^{-1}
$$
 for some matrix $T$ independent of $v\in \bbbz^k$.
Since $\rho_1(v)$ is a scalar matrix, the equation is equivalent to 
$$
 \rho_1(M\cdot v)\rho_2(\psi(v))=\rho_1(v),
$$
 i.e.
$$
 \rho_1((1-M)v)=\rho_2(\psi(v)).
$$
Note that $\hat \bbbz^k$ is isomorphic to the torus $T^k$,
 and the transformation $\rho_1\rightarrow \rho_1\circ(1-M)$
 is given by the action of the matrix $1-M$ on the torus $T^k$.
Therefore  the number of $\rho_1$ satisfying the last equation is the
degree
 of the map $(1-M)$ on the torus, i.e. $\mid \det(1-M)\mid $.
From this it follows that 
$$
 \#\fix(\hat\phi)
 =
 \#\fix(\hat\phi^{finite})\times \mid \det(1-M)\mid.
$$
As in the proof of Proposition 1
 we have $R(\phi ^{\infty})= \mid \det(1-M)\mid$.
Since $\phi^\infty$ is an endomorphism of an Abelian group we have 
 $\#\fix(\hat\phi^\infty)=R(\phi ^\infty)$.
Therefore
$$
 \#\fix(\hat\phi)
 =
 \#\fix(\hat\phi^{finite})\times \#\fix(\hat\phi^\infty). 
$$
\square
\medskip

As a consequence we have the following:

\begin{theorem}
If $\phi$ be any endomorphism of $G$ where $G$ is the direct product of a 
 finite group $F$ with a finitely generated free Abelian group, then
$$
 R(\phi)=\#\fix(\hat\phi)
$$
\end{theorem}

\proof
Since $\phi^{finite}$ is an endomorphism of a finite group, by 
Theorem 2 we have $R(\phi^{finite})=\#\fix(\hat\phi^{finite})$.
Since $\phi^\infty$ is an endomorphism of the finitely generated free 
Abelian group we have $R(\phi^\infty)=\#\fix(\hat\phi^\infty)$
 (see Theorem 2 in \cite{fh2}).
It follows from Propositions 1 and 2 that $R(\phi)=\#\fix(\hat\phi)$.
\square
\medskip

Let $W$ be the vector space of complex valued class functions on $F$.
The map $\phi$ induces a linear map $B:W\rightarrow W$
 defined as above in Theorem \ref{finitetrace}.

\begin{theorem}
\label{trace}
If $G$ is the direct product of a free Abelian and a finite group
 and $\phi$ an endomorphism of $G$.
Then the following trace formula holds:
\begin{equation}
 R(\phi)
 =
 (-1)^{r+p} \sum_{i=0}^k (-1)^i \tr (\wedge^i\phi^\infty\otimes B),
\end{equation}
 where $k$ is the rank of the free Abelian group $G/F$;
 $p$ is the number of $\mu\in\spec\phi^\infty$ such that $\mu <-1$,
 and $r$ the number of real eigenvalues of $\phi^\infty$ whose absolute
value is $>1$.
\end{theorem}

\proof
This follows from Theorems \ref{determinant} and \ref{finitetrace},
 Proposition \ref{product} and the formula
$$
 \tr (\wedge^i\phi^\infty)\cdot \tr (B)
 =
 \tr (\wedge^i\phi^\infty\otimes B).
$$
\square
\medskip

Finally let $\Gamma$ be a finitely generated, torsion free, nilpotent
group.
It is known \cite{mal} that any such group is a uniform, 
 discrete subgroup of some simply connected nilpotent Lie group $G$
 (uniform means that the coset space $G/ \Gamma$ is compact).
The space $M=G/ \Gamma$ is called a nil-manifold. 
Since $\Gamma=\pi_1(M)$ and $M$ is a $K(\Gamma,1)$,
 every endomorphism $\phi:\Gamma \to \Gamma$
 may be realized by a self-map $f:M\to M$ such that $f_*=\phi$
 and thus $R(f)=R(\phi)$.
Any endomorphism 
 $\phi:\Gamma \to \Gamma$ may be uniquely extended to an 
 endomorphism  $F: G\to G$.
Let $\tilde F:\tilde G\to \tilde G$ be the corresponding
 Lie algebra endomorphism induced from $F$.  

\begin{theorem}
\label{niltrace}
If $\phi$ an endomorphism of a finitely generated torsion free nilpotent 
 group $\Gamma$ then
\begin{equation}
 R(\phi)
 =
 (-1)^{r+p} \sum_{i=0}^m (-1)^i \tr \wedge^i\tilde F,
\end{equation}
 where $m$ is $\rank\Gamma = \dim M$;
 $p$ is the number of eigenvalues $\mu$ of $\tilde F$
 such that $\mu <-1$,
 and $r$ the number of real eigenvalues of $\tilde F$
 whose absolute value is $>1$.
\end{theorem} 

\proof 
Let $f:M\to M$ be a map realizing $\phi$ on a compact nil-manifold $M$ of
 dimension $m$.
We are assuming throughout this article
 that the Reidemeister number $R(f)=R(\phi)$ is finite. 
The finiteness of $R(f)$
 implies the non-vanishing of the Lefschetz number $L(f)$ \cite{fhw}. 
A strengthened version of Anosov's Theorem \cite{a}
 is proven in \cite{n} which states in particular,
 that if $L(f)\ne 0$ then $N(f)=|L(f)|=R(f)$. 
However it is known that $L(f)=\det (\tilde F -1)$ \cite{a}.
From this we have
$$ 
 R(\phi)
 =
 R(f)
 =
 |L(f)|
 =
 |\det (1- \tilde F)|=(-1)^{r+p}\det (1- \tilde F)
$$
$$
 =
 (-1)^{r+p} \sum_{i=0}^m (-1)^i \tr \wedge^i\tilde F.   
$$ 
\square
\medskip

\begin{theorem}(\cite{fhw})
Let $\Gamma$ be a finitely generated torsion free nilpotent group of 
rank $k$.
For any endomorphism $\phi:\Gamma \to \Gamma$ such that $R(\phi)$ is 
finite, there
exists an endomorphism $\psi:{\bbbz }^k\to {\bbbz}^k$ such that 
for all $n\in\N$,
$R(\phi^{n})=\#\fix \hat{\psi}^{n}$.
\end{theorem}

The following lemma is also useful for calculating Reidemeister numbers;
 it often allows one to reduce to the case that $\phi$ is an 
 isomorphism.

\begin{lemma}
Let $\phi:G\to G$ be any endomorphism
 of any group $G$, and let $H$ be a subgroup
 of $G$ with the properties
$$
  \phi(H) \subset H
$$
$$
  \forall x\in G \; \exists n\in \N \hbox{ such that } \phi^n(x)\in H.
$$
Then
$$
 R(\phi) = R(\phi_H),
$$
 where $\phi_H:H\to H$ is the restriction of $\phi$ to $H$.
\end{lemma}

{\it Proof .}
Let $x\in G$. Then there is an $n$ such that $\phi^n(x)\in H$. It is 
known that $x$ is $\phi$-conjugate
 to $\phi^n(x)$(see for example \cite{j}.
This means that the $\phi$-conjugacy class $\{x\}_\phi$
 of $x$ has non-empty intersection with $H$.

Now suppose that $x,y\in H$ are $\phi$-conjugate,
 ie. there is a $g\in G$ such that
 $$gx=y\phi(g).$$
We shall show that $x$ and $y$ are $\phi_H$-conjugate,
 ie. we can find a $g\in H$ with the above property.
First let $n$ be large enough that $\phi^n(g)\in H$.
Then applying $\phi^n$ to the above equation we obtain
 $$ \phi^n(g) \phi^n(x) = \phi^n(y) \phi^{n+1}(g). $$
This shows that $\phi^n(x)$ and $\phi^n(y)$ are $\phi_H$-conjugate.
On the other hand, one knows  that $x$ and $\phi^n(x)$ are
 $\phi_H$-conjugate, and $y$ and $\phi^n(y)$ are $\phi_H$ conjugate,
 so $x$ and $y$ must be $\phi_H$-conjugate.

We have shown that the intersection with $H$ of a 
 $\phi$-conjugacy class in $G$ is a $\phi_H$-conjugacy class
 in $H$.
We therefore have a map
$$
\begin{array}{cccc}
 Rest : & \cR(\phi) & \to & \cR(\phi_H)\\
        & \{x\}_\phi & \mapsto & \{x\}_\phi \cap H
\end{array}
$$
This clearly has the two-sided inverse
$$
 \{x\}_{\phi_H} \mapsto \{x\}_\phi.
$$
Therefore $Rest$ is a bijection and $R(\phi)=R(\phi_H)$.
\square

\begin{corollary}
 Let $H=\phi^n(G)$. Then $R(\phi) = R(\phi_H)$.
\end{corollary}

As an application of this we prove congruences for Reidemeister
 numbers.
Let $\mu(d)$, $d\in\bbbn$ be the M\"obius function,
 i.e.
$$
 \mu(d)
 =
 \left\{
 \begin{array}{ll}
 1 & {\rm if}\ d=1,  \\
 (-1)^k & \hbox{if $d$ is a product of $k$ distinct primes,}\\
 0 & \hbox{if $d$ is not square-free.} 
 \end{array}
 \right.
$$

\begin{theorem}
Let $\phi:G$ $\to G$ be an endomorphism of the group $G$ such that 
 all numbers $R(\phi^n)$ are finite, and let $H$ be a subgroup
 of $G$ with the properties
$$
 \phi(H) \subset H
$$
$$
 \forall x\in G \; \exists n\in \N \hbox{ such that } \phi^n(x)\in H.
$$
If one of the following conditions is satisfied:\\
(I) $H$ is a direct product of a finite group and a finitely generated 
free Abelian group,\\
 or\\
(II) $H$  is finitely generated, nilpotent and torsion free ,
then one has for all natural numbers $n$,
$$
 \sum_{d\mid n} \mu(d)\cdot R(\phi^{n/d}) \equiv 0 \mod n.
$$
\end{theorem}

\proof
From Theorem 4 and Lemma 2 it follows immediately that, in 
 case I, for every $n$ 
$$
 R(\phi^n)=\#\fix\left[ \hat{\phi_H}^n:\hat{H}\to\hat{H} \right].
$$

Let $P_n$ denote the number of periodic points of $\hat{\phi_H}$
 of least period $n$.
One sees immediately that
$$
 R(\phi^n)=\#\fix\left[ \hat{\phi_H}^n \right] = \sum_{d\mid n} P_d.
$$
Applying M\"obius' inversion formula, we have,
$$
 P_n = \sum_{d\mid n} \mu(d) R(\phi^{n/d}).
$$
On the other hand, we know that $P_n$ is always divisible
be $n$, because $P_n$ is exactly $n$ times the number of
$\hat{\phi_H}$-orbits in $\hat{H}$ of length $n$.
In the case II  when  $H$ is finitely generated, nilpotent and 
torsion free, we know from Theorem 7 that  there
exists an endomorphism $\psi:{\bf Z}^n\to {\bf Z}^n$ such that 
$R(\phi^n)=\#
\fix \hat{\psi^n}$.
 The proof then follows as in previous case.
\square

\section{ Reidemeister  and Nielsen zeta functions}  

In this section we reinterpret the results of the previous section
 in terms of the Reidemeister zeta function.
We show in the cases that we have considered that the zeta function is 
 a rational function with a functional equation.

% \begin{quotation} 
%  {\sc PROBLEM}. 
% For which groups and endomorphisms is the Reidemeister zeta function
% a rational function? When does it have a functional equation? 
% Is $R_\phi(z)$ an algebraic function? 
% \end{quotation}

\begin{theorem}
Let $G$ be the direct product of a free Abelian and a finite group
 and $\phi$ an endomorphism of $G$.
 Then $R_\phi(z)$ is a rational function and is equal to  
\begin{equation} 
 R_\phi(z) 
 = 
 \left( 
 \prod_{i=0}^k
 \det(1-\wedge^i\phi^\infty\otimes B \cdot\sigma\cdot z)^{(-1)^{i+1}} 
 \right)^{(-1)^r} 
\end{equation} 
 where matrix $B$ is as defined in Theorem 2, $\sigma=(-1)^p$, $p$ , $r$ 
 and $k$ are the constants described in Theorem 3. 
\end{theorem}
  
\proof
From Proposition 1 it follows that
 $R(\phi^n)=R((\phi^\infty)^n\cdot R((\phi^{finite})^n)$.  
From this formula,    Theorem 2  and  3 we have the trace 
 formula for $R(\phi^n)$: 
\begin{eqnarray*} 
 R(\phi^n) & = & (-1)^{r+pn} \sum_{i=0}^k (-1)^i \tr 
\wedge^i(\phi^\infty)^n \cdot \tr B^n \\ 
           & = &  (-1)^{r+pn} \sum_{i=0}^k (-1)^i \tr 
(\wedge^i(\phi^\infty)^n\otimes B^n) \\ 
         & = & (-1)^{r+pn} \sum_{i=0}^k (-1)^i \tr 
(\wedge^i\phi^\infty\otimes B )^n. 
\end{eqnarray*} 
  
We now calculate directly
\begin{eqnarray*} 
 R_\phi(z)
 & = &
 \exp\left(\sum_{n=1}^\infty \frac{R(\phi^n)}{n} z^n\right) \\ 
 & = &
 \exp\left(\sum_{n=1}^\infty  
                \frac{(-1)^{r} 
                \sum_{i=0}^k (-1)^i
                \tr (\wedge^i\phi^\infty\otimes B )^n}{n} (\sigma \cdot
z)^n\right) \\ 
 & = &
 \left(\prod_{i=0}^k\left( \exp\left(\sum_{n=1}^\infty  
                \frac{1}{n}\tr (\wedge^i\phi^\infty\otimes B )^n\cdot
(\sigma \cdot z)^n\right) 
                \right)^{(-1)^i} \right)^{(-1)^r} \\ 
 & = &
 \left(\prod_{i=0}^k
 \det\left(
 1-\wedge^i\phi^\infty\otimes B\cdot\sigma\cdot z\right)^{(-1)^{i+1}}
 \right)^{(-1)^r}. 
\end{eqnarray*} 
\square

\begin{theorem} 
If $\Gamma$ is a finitely generated torsion free nilpotent group and 
$\phi$ an
endomorphism of $\Gamma$ .Then $R_\phi(z)$ is a rational function and 
is equal to  
\begin{equation} 
 R_\phi(z) 
 = 
 \left( 
 \prod_{i=0}^m
 \det(1-\wedge^i\tilde F \cdot\sigma\cdot z)^{(-1)^{i+1}} 
 \right)^{(-1)^r}  
\end{equation} 
where $\sigma=(-1)^p$,$p$ , $r$, $m$ and  $\tilde F$ is defined in 
theorem 6 
\end{theorem} 
  
\proof
If we repeat the proof of the Theorem 4 for $\phi^n$ instead $\phi$ 
we obtain  
that $ R(\phi^n)=(-1)^{r+pn}\det(1-\tilde F)$( we suppose that 
Reidemeister numbers  $ R(\phi^n)$ are finite for all $n$).Last 
formula implies the trace formula for $R(\phi^n)$:
$$ 
R(\phi^n) = (-1)^{r+pn} \sum_{i=0}^m (-1)^i \tr (\wedge^i\tilde F)^n
$$ 
From this we have  formula 6 immediately by direct calculation as in 
Theorem 9.  
\square

\begin{corollary} 
Under the assumptions of Theorem 6,
 the zeros and poles of the Reidemeister zeta function
 are reciprocals of eigenvalues of one of the maps
$$
 \wedge^i \tilde F : \wedge^i \tilde G \longrightarrow \wedge^i \tilde G
 \qquad
 0\leq i\leq \rank\;\Gamma
$$ 
\end{corollary}

We now turn our attention to the functional equations satisfied by 
these zeta functions. These equations are consequences of a duality 
on the vector spaces with which we have obtained trace formulae.

\begin{theorem}
Let $\Gamma$ be a finitely generated, torsion free, nilpotent group
 and $\phi$ an endomorphism of $\Gamma$.
The Reidemeister zeta function $R_\phi(z)$ satisfies the
 following functional equation: 
\begin{equation} 
 R_{\phi}\left(\frac{1}{dz}\right) 
 = 
 \epsilon_1 \cdot R_\phi(z)^{(-1)^{\rank \Gamma}},
\end{equation} 
 where $d=\det\tilde F$ and $\epsilon_1$ is a non-zero complex constant. 
\end{theorem} 

\proof
Via the natural non-singular pairing
 $(\wedge^i \tilde F)\otimes(\wedge^{m-i} \tilde F) \rightarrow\bbbc$ 
 the operators $\wedge^{m-i}\tilde F$ and $d\cdot (\wedge^i \tilde F)^{-1}$
 are adjoint to each other.

Consider an eigenvalue $\lambda$ of $\wedge^i\tilde F$. 
By Theorem 10, This contributes a term
$$ 
 \left((1-\frac{\lambda\sigma}{dz})^{(-1)^{i+1}}\right)^{(-1)^r} 
$$ 
 to $R_\phi\left(\frac{1}{dz}\right)$. 
This may be rewritten as follows:
$$ 
 \left( 
 \left( 1 - \frac{d\sigma z}{\lambda} \right)^{(-1)^{i+1}}\cdot 
 \left( \frac{-dz}{\lambda\sigma} \right)^{(-1)^i} 
 \right)^{(-1)^r}.
$$ 
Note that $\frac{d}{\lambda}$ is an eigenvalue of $\wedge^{m-i}\tilde F$. 
Multiplying these terms together we obtain, 
$$ 
 R_\phi\left(\frac{1}{dz}\right) 
 = 
 \left( 
 \prod_{i=1}^m
 \prod_{\wedge^{(i)}\in\spec\wedge^i\tilde F}
 \left(\frac{1}{\wedge^{(i)}\sigma}\right)^{(-1)^i}
 \right)^{(-1)^r}
 \times
 R_\phi(z)^{(-1)^m}.
$$
The variable $z$ has disappeared because
$$
 \sum_{i=0}^m (-1)^i \dim\wedge^i\tilde G
 =
 \sum_{i=0}^m (-1)^i  {C_m }^i = 0.
$$
\square

\subsection[Continuous maps]{The Reidemeister zeta functions of a continuous
map.}

Now suppose $X$ is a compact, connected polyhedron
 with universal cover $\tilde X$
 and let $f:X\rightarrow X$ is a continuous self-map.
We choose a lifting $\tilde{f}:\tilde{X}\rightarrow\tilde{X}$ of $f$ as
reference.
Let $\Gamma$ be the group of
 covering translations of $\tilde{X}$ over $X$.
Then every lifting of $f$ can be written uniquely
 as $\gamma\circ \tilde{f}$, with $\gamma\in\Gamma$.
So elements of $\Gamma$ serve as coordinates of
 liftings with respect to the reference $\tilde{f}$.
Now for every $\gamma\in\Gamma$ the composition $\tilde{f}\circ\gamma$
 is a lifting of $f$ so there is a unique $\gamma^\prime\in\Gamma$
 such that $\gamma^\prime\circ\tilde{f}=\tilde{f}\circ\gamma$.
This correspondence $\gamma\mapsto\gamma^\prime$ is an 
 endomorphism of $\Gamma$.
As such it is dependent on the choice of $\tilde{f}$.
However the class of the endomorphism modulo inner endomorphisms of 
 $\Gamma$ depends only on the map $f$.

\begin{definition}
The endomorphism $\tilde{f}_*:\Gamma\rightarrow\Gamma$ determined
 by the lifting $\tilde{f}$ of $f$ is defined by
$$
 \tilde{f}_*(\gamma)\circ\tilde{f} = \tilde{f}\circ\gamma.
$$
\end{definition}

\begin{lemma}\cite{j}
Lifting classes of $f$ are in one to one correspondence with
$\tilde{f}_*$-conjugacy classes in $\pi$,
the lifting class $[\gamma\circ\tilde{f}]$ corresponding
to the $\tilde{f}_*$-conjugacy class of $\gamma$.
We therefore have $R(f) = R(\tilde{f}_*)$.
\end{lemma}

Using this lemma we may apply all previous theorems  to the 
Reidemeister zeta functions of continuous maps.

\subsection{Convergence of the Nielsen Zeta Function}
 
We denote by $R$ the radius of convergence of the Nielsen zeta 
 function $N_f(z)$,
 and by $h(f)$ the topological entropy of  a continuous map $f$.
Let $h$ be the infimum of $h(g)$ where $g$ ranges over maps
 with the same homotopy type as $f$.
The following is known.

\begin {theorem}[\cite{fp}]
With the above notation $R \geq \exp(-h) > 0$.
\end{theorem}

The radius of convergence $R$ is directly connected with 
 asymptotic Nielsen number which is defined to be the growth rate of 
  the sequence $\{N(f^{n})\}$ of   the Nielsen numbers of the iterates of $f$.
In this section we give a new proof that $R$ is positive,
 and give an exact algebraic lower bound for $R$ using the trace 
 formulae (13) and (14) for generalized Lefschetz numbers.

% \subsubsection[Reidemeister trace formula]{The Reidemeister trace 
% formula for generalized Lefschetz numbers}

We begin by recalling some known facts (see
\cite{re},\cite{w},\cite{fahu},\cite{j1}).
% We shall use this results to estimate the radius of convergence of 
%  the Nielsen zeta function. 
The fundamental group $\pi=\pi_1(X,x_0)$ splits into 
 $\tilde{f}_*$-conjugacy classes.
Let $\pi_f$ denote the set of $\tilde{f}_*$-conjugacy classes,
 and $\bbbz\pi_f$ denote the Abelian group freely generated by $\pi_f$.
We shall use the bracket notation 
 $a\mapsto [a]$ for both projections $\pi\to \pi_f$ and
$\bbbz\pi\to\bbbz\pi_f$. 
Let $x$ be a fixed point of $f$.
Choose a path $c$ from $x_0$ to $x$.
The $\tilde{f}_*$-conjugacy class in $\pi$ of the loop $c\cdot (f\circ
c)^{-1}$,
 which is evidently independent of the choice of $c$, 
 is called the coordinate of $x$.
Two fixed points are in the same fixed point class $F$
 iff they have the same coordinates.
This $\tilde{f}_*$-conjugacy class is thus called the coordinate of the 
 fixed point class $F$ and denoted $cd_{\pi}(F,f)$.
The generalized Lefschetz number or the Reidemeister trace \cite{re}
 is defined as
\begin{equation} 
 L_{\pi}(f)
 :=
 \sum_{F}\ind (F,f)\cdot cd_{\pi}(F,f)
 \in
 \bbbz\pi_f , 
\end{equation} 
 the summation being over all essential fixed point classes $F$ of $f$.
The Nielsen number $N(f)$ is the number of non-zero terms in $L_{\pi}(f)$,
 and the indices of the essential fixed point classes
 appear as the coefficients in $L_{\pi}(f)$.
This invariant used to be called the Reidemeister trace  
 because it can be computed as an alternating sum of traces on the chain
level   
 as follows (\cite{re},\cite {w}). 
Assume that $X$ is a finite cell complex and  $f:X\to X$ is a cellular map.

A cellular decomposition $\{e_j^d\}$ of $X$ lifts to a $\pi$-invariant
cellular structure
 on the universal covering $\tilde X$.
Choose an arbitrary lift $\tilde{e}_j^d$ for each $e_j^d$.
These lifts constitute a free $\bbbz\pi$-basis
 for the cellular chain complex of $\tilde{X}$. 
The lift $\tilde{f}$ of $f$ is also a cellular map.
In every dimension $d$, the cellular chain map $\tilde{f}$  
 gives rise to a $\bbbz\pi$-matrix $\tilde{F}_d$ with respect to 
 the above basis, i.e. $\tilde{F}_d=(a_{ij})$
 if $\tilde{f}(\tilde{e}_i^d)=\sum_{j}a_{ij}\tilde{e}_j^d $,
 where $a_{ij}\in \bbbz\pi$.
Then we have the Reidemeister trace formula
\begin{equation} 
 L_{\pi}(f)
 =
 \sum_{d}(-1)^d[\Tr \tilde{F}_d]
 \in
 \bbbz\pi_f. 
\end{equation}

We now describe alternative approach to the Reidemeister trace 
 formula proposed by Jiang \cite{j1}.
This approach is useful when we study the periodic points of $f$.  
  The mapping torus $T_f$ of $f:X\rightarrow X$ is the space obtained 
 from $X\times [0,\infty )$ by identifying $(x,s+1)$ with $(f(x),s)$ 
 for all $x\in X,s\in [0 ,\infty )$.
On $T_f$ there is a natural semi-flow
 $\phi :T_f\times [0,\infty )\rightarrow T_f,\phi_t(x,s)=(x,s+t)$
 for all $t\geq 0$.
The function $f:X\rightarrow X$ is the return map of the semi-flow $\phi$.
A point $x\in X$ and a positive number $\tau$ determine the orbit curve
 $\phi _{(x,\tau)}:={\phi_t(x)}_{0\leq t \leq \tau}$ in $T_f$. 
Take the base point $x_0$ of $X$ as the base point of $T_f$.
It is known that the fundamental group $H:=\pi_ 1(T_f,x_0)$
 is obtained from $\pi$ by adding a new generator $z$
 and adding the relations $z^{-1}gz=\tilde f_*(g)$
 for all $g\in \pi =\pi _1(X,x_0)$.
Let  $H_c$ denote the set of conjugacy classes in $H$
 and let $\bbbz H$ be the integral group ring of $H$.
Define also $\bbbz H_c$ to be the free Abelian group with basis $H_c$.
We again use the bracket notation $a\rightarrow [a]$ for both projections
 $H \rightarrow H _c $ and $ \bbbz H \rightarrow \bbbz H_c $.
If $F^n$ is a fixed point class of $f^n$, then  
 $f(F^n)$ is also fixed point class of $f^n$
 and $\ind (f(F^n),f^n)=\ind (F^n,f^n)$. 
Thus $f$ acts as an index-preserving permutation
 among fixed point classes of $f^n$.
We define an $n$-orbit class $O^n$ of $f$ to 
 be the union of elements of an orbit of this action.
In other words, two points $x,x'\in \Fix (f^n)$ are said to be in 
 the same $n$-orbit class of $f$ if and only if some $f^i(x)$ and some 
 $f^j(x')$ are in the same fixed point class of $f^n$.
The set $\Fix (f^n)$ is a disjoint union of $n$-orbit classes.
A point $x$ is a fixed point of $f^n$   
% or a periodic point of period $n$
 if and only if $\phi_{(x,n)}$ is a closed curve.  
The free homotopy class of the closed curve $\phi_{(x,n)}$ will be 
 called the $H$-coordinate of point $x$,
 written $cd_{H }(x,n)=[\phi _{(x,n)}]\in H_c$.
It follows that periodic points $x$  
 of period $n$ and $x'$ of period $n'$ have the same $H$-coordinate 
 if and only if $n=n'$, and $x$ and $x'$ belong to the same $n$-orbits class of $f$. 

% Thus it is possible equivalently define $x,x'\in \Fix (f^n)$ to be 
%  in the same $n$-orbit class if and only if they have the same 
%  $H-$coordinate.  
Recently, Jiang \cite{j1} has considered generalized Lefschetz 
 number with respect to $H$:
\begin{equation} 
 L_{H}(f^n)
 :=
 \sum_{O^n} \ind(O^n,f^n) \cdot cd_{H }(O^n)
 \in
 \bbbz H_c.
\end{equation} 
He proved the following trace formula: 
\begin{equation} 
 L_{H}(f^n)
 =
 \sum_{d}(-1)^d[\Tr (z\tilde{F}_d)^n]
 \in
 \bbbz H_c, 
\end{equation} 
 where $\tilde{F}_d$ are the $\bbbz \pi$-matrices defined in (9) and 
 $z\tilde{F}_d$ is regarded as a $\bbbz H$-matrix.

For any set $S$ let $\bbbz S$ denote the free Abelian group with basis $S$.
Define a norm on $\bbbz S$ by
\begin{equation}
 \|\sum_i k_is_i\|
 :=
 \sum_i \mid k_i\mid \in \bbbz,
\end{equation}
 where the $s_i$ in $S$ are all different.

For a $\bbbz H$-matrix $A=(a_{ij})$,
 define its norm by
 $\| A\| := \sum_{i,j}\| a_{ij} \|$.
Then we have inequalities
 $\| AB \|\leq \| A\ |\cdot\| B \|$
 when $A,B$ can be multiplied,
 and $\| \tr A \|\leq \| A \|$ when $A$ is a square matrix.
For a matrix  $A=(a_{ij})$ in $\bbbz S$,
 its matrix of norms is defined to be the 
 matrix $A^{norm}:=(\| a_{ij} \|)$ which is a matrix of non-negative 
 integers.
In what follows, the set $S$ will be $\pi$, $H$ or
 $H_c$.
We denote by $s(A)$ the spectral radius of $A$,
 $s(A)=\lim_n (\| A^n \| |)^{\frac{1}{n}}$,
 which coincides with the largest modulus of an eigenvalue of $A$.

\begin{theorem}
For any continuous map $f$ of any compact polyhedron $X$ into 
 itself the Nielsen zeta function has positive radius of convergence 
 $R$, which admits following estimations 
\begin{equation}
 R\geq \frac{1}{\max_d \| z\tilde F_d \|} > 0
\end{equation}
 and
\begin{equation}
 R\geq \frac{1}{\max_d s(\tilde F_d^{norm})} > 0, 
\end{equation}
 where $\tilde F_d$ is as in (9).
\end{theorem}

\proof
By the homotopy type invariance of the invariants we can suppose that 
 $f$ is a cellular map of a finite cell complex.
By definition the Nielsen number $N(f^n)$ is the number of non-zero terms
in 
 $L_{\pi}(f^n)$.
The norm $\| L_{H}(f^n) \|$ is the sum of absolute 
 values of the indices of all the $n$-orbits classes $O^n$.
It equals $\| L_{\pi}(f^n) \|$, the sum of absolute values of the indices
of 
 all the fixed point classes of $f^n$, because any two fixed point 
 classes of $f^n$ contained in the same $n$-orbit class $O^n$ must 
 have the same index.
From this we have 
 $N(f^n)
 \leq\|L_{\pi}(f^n)\|
 =\|L_H(f^n)\|
 =\|\sum_d(-1)^d[\tr(z\tilde F_d)^n]\|
 \leq\sum_d\|[\tr(z\tilde F_d)^n]\|
 \leq\sum_d\|\tr(z\tilde F_d)^n\|
 \leq\sum_d\|(z\tilde F_d)^n\|
 \leq\sum_d\|(z\tilde F_d)\|^n$
 (see \cite{j1}).
The radius of convergence $R$ is given by the Cauchy-Hadamard formula:
$$
 \frac{1}{R}
 =
 \limsup_n (\frac{N(f^n)}{n})^{\frac{1}{n}}
 =
 \limsup_n (N(f^n))^{\frac{1}{n}}.
$$
Therefore we have:
$$
 R
 =
 \frac{1}{\limsup_n (N(f^n))^{\frac{1}{n}}}
 \geq
 \frac{1}{\max_d \|z\tilde F_d \|}
 >
 0.
$$
The inequalities:
$$
 N(f^n)\leq\|L_{\pi}(f^n)\|\!=\!\|L_H(f^n)\|
 \!=\!
 \|\sum_d(-1)^d[\tr(z\tilde F_d)^n]\|\!\leq\sum_d\| [\tr(z\tilde
F_d)^n]\|\!
$$
$$
 \leq
 \sum_d{||}\tr(z\tilde F_d)^n{||}
 \leq
 \sum_d\tr((z\tilde F_d)^n)^{norm}
 \leq
 \sum_d\tr((z\tilde F_d)^{norm})^n
$$
$$
 \leq
 \sum_d\tr((\tilde F_d)^{norm})^n
$$
 together with the definition of spectral radius give the bound:
$$
 R
 =
 \frac{1}{\limsup_n (N(f^n))^{\frac{1}{n}}}
 \geq
 \frac{1}{\max_d s(\tilde F_d^{norm})}
 >
 0.
$$
\square

\begin{example}
Let $X$ be a surface with boundary, and  $f:X\rightarrow X$ a 
 continuous map.
Fadell and Husseini \cite{fahu} devised a method of computing the 
 matrices of the lifted chain map for surface maps.
Suppose $\{a_1, \ldots ,a_r\}$ is a free basis for $\pi_1(X)$.
Then $X$ has the homotopy type of a bouquet $B$ of $r$ circles which can be
decomposed into one 0-cell and $r$ 1-cells corresponding to the 
 $a_i$, and $f$ has the homotopy type of a cellular map  
 $g:B\rightarrow B$.
By the homotopy invariance, we may replace $f$ 
 by $g$ in computations.
The homomorphism $\tilde f_*:\pi_1(X)\rightarrow \pi_1(X)$
 induced by $f$ or $g$ is determined by the images
 $b_i=\tilde f_*(a_i)$, $i=1,\ldots,r $.
The fundamental group $\pi_1(T_f)$ has a presentation 
 $\pi_1(T_f)=<a_1,...,a_r,z| a_iz=zb_i, i=1,..,r>$.
Let
$$
 D=(\frac{\partial b_i}{\partial a_j})
$$
 be the Jacobian in Fox calculus (see \cite{bi}).
Then, as pointed out in \cite{fahu},
 the matrices of the lifted chain map $\tilde g$ are
$$
 \tilde F_0=(1),\qquad
 \tilde F_1=D=(\frac{\partial b_i}{\partial a_j}).
$$
Now, we may find bounds for the radius $R$ using the above theorem.
\end{example}

\section{Reidemeister Torsion as a Special Value}

Reidemeister torsion is an algebraically defined quantity associated
 to an acyclic cochain complex.
It is used to distinguish between complexes with the same homology.
Roughly speaking, if Euler characteristic is regarded as a
 graded version generalization of dimension,
 then Reidemeister torsion may be viewed as a graded version of
 the absolute value of the determinant.
More precisely, let $d^i:C^i\rightarrow C^{i+1}$ be a bounded,
 acyclic cochain complex of finite dimensional complex vector spaces.
There is a chain contraction $\delta^i:C^i\rightarrow C^{i-1}$
 ie. a linear map such that $d\circ\delta + \delta\circ d = id$.
We have linear maps
 $ (d+\delta)_+ : C^+:=\oplus C^{2i} \rightarrow C^- :=\oplus C^{2i+1}$
 and $(d+\delta)_- : C^- \rightarrow C^+$.
Since $(d+\delta)^2 = id + \delta^2$ is unipotent,
 it follows that $(d+\delta)_+$ is bijective.
Given positive densities $\Delta_i$ on $C^i$,
 we may define with respect to these densities,
 $\tau(C^*,\Delta_i):= \mid \det(d+\delta)_+\mid$ (see \cite{fri2}).

Reidemeister torsion is used in the following geometric setting.
Let $K$ is a finite simplicial complex and $E$ a flat, finite dimensional,
 complex vector bundle over $K$ with fibre $V$.
Denote by $\rho_{E}:\pi_{1}(K)\to \GL(V)$ the holonomy of $E$.
Suppose now that one has on each fibre of $E$ a positive density
 which is locally constant on $K$.
In terms of $\rho_E$ this assumption means $\mid\det\rho_E\mid=1$.
The cochain complex $C^i(K;E)$ with coefficients in $E$
 may be identified with the direct sum of copies
 of $V$ associated to each $i$-cell $\sigma$ of $K$.
This identification is achieved by choosing a base point in each
 component of $K$ and a base point in each $i$-cell.
By choosing a flat density on $E$ we obtain a
 preferred density $\Delta_i$ on each $C^i(K,E)$.
The Reidemeister torsion $\tau(K;E)$ is then defined
 to be the positive real number $\tau(C^*(K;E),\Delta_i)$.
This has many nice properties.
It is invariant under subdivisions of $K$.
Thus for a smooth manifold,
 one may unambiguously define $\tau(K;E,D_i)$
 to be the torsion of any smooth triangulation of $K$.
However, Reidemeister torsion
 is not an invariant under a general homotopy equivalence.
This was in fact the reason for its introduction.

In the case that $K$ is the circle $S^{1}$,
 let $A$ be the holonomy of a generator of its fundamental group 
 $\pi_1(S^1)$.
Then $E$ is acyclic if and only if $I-A$ is invertible and then
\begin{equation}
 \tau(S^1;E)
 =
 \mid \det(I-A) \mid                                     
\end{equation}
Note that the choice of generator is irrelevant as
 $I-A^{-1}=(-A^{-1})(I-A)$ and $|det(-A^{-1})|=1$.

It might be expected that the Reidemeister torsion counts something 
 geometric (like the Euler characteristic).
D. Fried showed that it counts the periodic orbits of a flow and the 
 periodic points of a map.
We will show that Reidemeister torsion counts the periodic point classes
 of a map (fixed point classes of the iterations of the map).

Some further properties of $\tau$ describe its behavior under bundles. 

Suppose $p: X\rightarrow B$ is a simplicial bundle with fiber $F$,
 where  $F, B, X$ are all finite complexes
 and $p^{-1}$ lifts subcomplexes of $B$ to subcomplexes of $X$. 
We assume here that $E$ is a flat, complex vector bundle over $B$ . 
We form its pullback $p^*E$ over $X$.
Note that the vector spaces $H^i(p^{-1}(b),\bbbc)$ with
 $b\in B$ form a flat vector bundle over $B$,
 which we denote $H^i F$.
The integral lattice in $H^i(p^{-1}(b),\bbbr)$
 determines a flat density by the condition
 that the covolume of the lattice is $1$.
Suppose that the bundle $E\otimes H^i F$ is acyclic for all $i$.
Under these conditions D. Fried \cite{fri2} has shown that the bundle
 $p^* E$ is acyclic, and
\begin{equation}
 \tau(X;p^* E)
 =
 \prod_i \tau(B;E\otimes H^i F)^{(-1)^i}.
\end{equation}
\medskip

Now let $f:X\rightarrow X$ be a homeomorphism of
 a compact polyhedron $X$,
 and let $T_f := (X\times I)/(x,0)\sim(f(x),1)$ be the mapping torus of
$f$.
We shall consider the bundle $p:T_f\rightarrow S^1$ over the circle $S^1$.
Given a flat, complex vector bundle $E$ with 
 finite dimensional fibre over the base $S^1$,
 we form its pullback $p^*E$ over $T_f$.
The vector spaces $H^i(p^{-1}(b),\bbbc)$ with
 $b\in S^1$ also form flat vector bundles over $S^1$,
 which we denote $H^i F$.
The integral lattice in $H^i(p^{-1}(b),\bbbr)$ determines a flat density
 by the condition that the covolume of the lattice is $1$.
Suppose that the bundle $E\otimes H^i F$ is acyclic for all $i$.
Under these conditions, D. Fried \cite{fri2} has shown that
 $p^* E$ is also acyclic and the following holds:
\begin{equation}
 \tau(T_f;p^* E)
 =
 \prod_i \tau(S^1;E\otimes H^i F)^{(-1)^i}.
\end{equation}
Let $g$ be the preferred generator of the group $\pi_1 (S^1)$
 and let $A=\rho(g)$, where $\rho:\pi_1 (S^1)\rightarrow GL(V)$.
The holonomy around $g$ of the bundle $E\otimes H^i F$
 is then $A\otimes f^*_i$.
 
Since $\tau(E)=\mid\det(I-A)\mid$ it follows from (24) that
\begin{equation}
 \tau(T_f;p^* E)
 =
 \prod_i \mid\det(I-A\otimes f^*_i)\mid^{(-1)^i}.
\end{equation}
Now consider the special case in which $E$ is one-dimensional,
 so the holonomy $A$ is a complex scalar $\lambda$ with absolute value one.
In terms of the rational function $L_f(z)$ we have \cite{fri2}:
\begin{equation}
 \tau(T_f;p^* E)
 =
 \prod_i \mid\det(I-\lambda .f^*_i)\mid^{(-1)^i}
 =
 \mid L_f(\lambda)\mid^{-1}
\end{equation}

\begin{theorem}
Let $\phi:G\to G$ be a group automorphism, where  $G$ is the direct 
 sum of a finite group with a finitely generated Abelian group, 
 then
$$
 \tau\left(T_{\hat{\phi}};p^*E\right)
 =
 \mid L_{\hat{\phi}}(\lambda) \mid^{-1}
 =
 \mid R_{\phi}(\sigma\lambda) \mid^{(-1)^{r+1}},
$$
 where $\lambda$ is the holonomy of $E$ around $S^1$
 and $r$ and $\sigma$ are the  constants described in Theorem 3.
\end{theorem}
 
\proof
We  know from Theorem 4
 that $R(\phi^n)$ is the number of fixed points
 of the map $\hat{\phi}^{n}$.
It remains to show that
 the number of fixed points of $\hat{\phi}^n$ is equal to the
 absolute value of its Lefschetz number.
We assume without loss of generality that $n=1$.
We are assuming throughout that $R(\phi)$ is finite, so
 the fixed points of $\hat{\phi}$ form a discrete set.
We therefore have
$$
 L(\hat{\phi})
 =
 \sum_{x\in\fix\hat{\phi}} \ind(\hat{\phi},x).
$$
Since $\phi$ is a group endomorphism, the trivial representation 
 $x_0\in\hat{G}$ is always fixed.
Let $x$ be any fixed point of $\hat{\phi}$.
Since $\hat{G}$ is union of tori $\hat{G}_0,...,\hat{G}_t$
 and $\hat{\phi}$ is a linear map, we can
 shift any two fixed points onto one another without altering the map 
 $\hat\phi$.
Therefore all fixed points have the same index.
It is now sufficient to show that $\ind(\hat{\phi},x_0)=\pm 1$.
This follows because the map on the torus
$$
 \hat{\phi}:\hat{G}_0\to\hat{G}_0
$$
 lifts to a linear map of the universal cover,
 which is an euclidean space.
The index is then the sign of the determinant
 of the identity map minus this lifted map.
This determinant cannot be zero, because $1-\hat{\phi}$ must have finite
 kernel by our assumption that the Reidemeister number of $\phi$ is
 finite (if $\det(1-\hat{\phi})=0$ then the kernel of $1-\hat{\phi}$
 is a positive dimensional subspace of $\hat{G}$, and therefore infinite).
\square

\section{Concluding remarks and open questions}
 
For the case of almost nilpotent groups (ie. groups with polynomial growth,
 in view of Gromov's Theorem \cite{gromov}) we believe that some 
 power of the Reidemeister zeta function is a rational function
 and that the congruences for the Reidemeister numbers are also true.
We intend to prove this conjectures by
 identifying the Reidemeister number on the nilpotent part of the group
 with the number of fixed points in the direct sums of the duals of the
 quotients of successive terms in the central series.
We then hope to show that the Reidemeister number 
 of the whole endomorphism is a sum of numbers of orbits of such fixed 
 points under the action of the finite quotient group
 (ie the quotient of the whole group by the nilpotent part).
The situation for groups with exponential growth is very different.
There one can expect the Reidemeister number to be infinite as long 
 as the endomorphism is injective.

\bigskip
 
 Institut f\"ur  Mathematik,E.-M.-Arndt- Universit\"at  Greifswald 
 
Jahn-strasse 15a, D-17489 Greifswald, Germany.
 
{\it E-mail address}: felshtyn@rz.uni-greifswald.de

\bigskip

Department of Mathematics, University College London,

London WC1E 6BT

{\it E-mail address}: rih@math.ucl.ac.uk

\end{document}